\documentclass[a4paper]{article}
\usepackage{amsmath,amssymb,amsfonts}
\usepackage[all]{xy}
\def\mar[#1]{\ar@{-}[#1]|-{\object@{<}}}

\CompileMatrices

\def\un{{\bf 1}}
\def\zero{\{0\}}
\def\mpoint{\;.}
\def\mvirg{\;,}

\def\mpn{\medskip\par\noindent}

\def\mmpn{\vskip 1em minus 1em\par\noindent}

\def\normal{\trianglelefteq}

\def\medp{\medskip\par}
\def\smp{\smallskip\par}
\def\Res{\hbox{\rm Res}}

\def\Hom{\hbox{\rm Hom}}

\def\Aut{\hbox{\rm Aut}}

\def\Out{\hbox{\rm Out}}
\def\Inf{\hbox{\rm Inf}}
\def\Def{\hbox{\rm Def}}

\def\Id{\hbox{\rm Id}}
\def\Im{\hbox{\rm Im}}
\def\Ker{\hbox{\rm Ker}}

\def\Defres{\hbox{\rm Defres}}

\def\Indinf{\hbox{\rm Indinf}}

\def\dom{\backslash}
\newcommand{\flh}[2]{\mathop{\hbox to  4ex{\rightarrowfill}}_{#2}^{#1}\limits}

\newcommand{\sumb}[2]{\sum_{{\scriptstyle #1}\atop {\scriptstyle #2}}}
\newcommand{\oplusb}[2]{\mathop{\bigoplus}_{{\scriptstyle  #1}\atop{\scriptstyle #2}}}

\newcommand{\ressort}[1]{\hskip #1 plus #1 minus #1}
\def\pf{\noindent{\bf Proof. }}
\def\findemo{~\leaders\hbox to 1em{\hss\  \hss}\hfill~\raisebox{.5ex}{\framebox[1ex]{}}\smp}
\newcommand{\carre}[8]{\begin{array}{ccc}
#1&\mathop{\hbox to  12mm{\rightarrowfill}}^{\displaystyle{#2}}\limits&#3\\
\llap{$\displaystyle{#4}$}\left\downarrow\vbox to 6mm{}\right. & &  \left\downarrow\vbox to 6mm{}\right.\rlap{$\displaystyle{#5}$}\\
\ & &\\#6&\mathop{\hbox to 12mm{\rightarrowfill}}_{\displaystyle  #7}\limits&#8\\
\end{array}}
\newcommand{\limind}[1]{\mathop{\lim}_{\displaystyle\longrightarrow\atop  \scriptstyle{#1}}\limits}
\newcommand{\limproj}[1]{\mathop{\lim}_{\displaystyle\longleftarrow\atop  \scriptstyle{#1}}\limits}

\makeatletter
\renewenvironment{enumerate}{\ifnum \@enumdepth >3 \@toodeep\else
      \advance\@enumdepth \@ne
      \edef\@enumctr{enum\romannumeral\the\@enumdepth}\list
      {\csname  label\@enumctr\endcsname}{\setlength{\topsep}{1ex}\setlength{\itemsep}{0 pt}\usecounter
        {\@enumctr}\def\makelabel##1{\hss\llap{##1}}}\fi}{\endlist}
\renewenvironment{itemize}{\ifnum \@itemdepth >3 \@toodeep\else  \advance\@itemdepth \@ne
\edef\@itemitem{labelitem\romannumeral\the\@itemdepth}%
\list{\csname\@itemitem\endcsname}{\setlength{\topsep}{1ex}\setlength{\itemsep}{0pt}\def\makelabel##1{\hss\llap{##1}}}\fi}
{\endlist}

\def\@sect#1#2#3#4#5#6[#7]#8{\ifnum #2>\c@secnumdepth
     \let\@svsec\@empty\else
     \refstepcounter{#1}\edef\@svsec{\csname the#1\endcsname .\hskip  .5em}\fi
     \@tempskipa #5\relax
      \ifdim \@tempskipa>\z@
        \begingroup #6\relax
          \@hangfrom{\hskip #3\relax\@svsec}{\interlinepenalty \@M  #8\par}%
        \endgroup
       \csname #1mark\endcsname{#7}\addcontentsline
         {toc}{#1}{\ifnum #2>\c@secnumdepth \else
                      \protect\numberline{\csname the#1\endcsname}\fi
                    #7}\else
        \def\@svsechd{#6\hskip #3\relax  
                   \@svsec #8\csname #1mark\endcsname
                      {#7}\addcontentsline
                           {toc}{#1}{\ifnum #2>\c@secnumdepth \else
                             \protect\numberline{\csname  the#1\endcsname}\fi
                       #7}}\fi
     \@xsect{#5}}

\def\section{\pagebreak[3]\setcounter{prop}{0}\setcounter{equation}{0}\@startsection{section}{1}{\z@}{6ex plus  9ex}{3ex}{\center\reset@font\large\bf}}
\def\subsubsection{\@startsection{subsubsection}{3}{\z@}{4ex plus  6ex}{-1em}{\reset@font\it}}
\makeatother
\newcommand{\masubsect}[1]{\medskip\par\noindent\pagebreak[3]\refstepcounter{subsection}\refstepcounter{prop}{\bf \thesection.\arabic{prop}.\ \ #1.\ }}
\def\Z{\mathbb{Z}}

\def\F{\mathbb{F}}

\def\theprop{\thesection.\arabic{prop}}

\renewenvironment{equation}{\refstepcounter{subsection}\refstepcounter{prop}$$}{\leqno{\bf (\theprop)}$$}

\newenvironment{rem}[1]{\refstepcounter{subsection}\refstepcounter{prop} \mpn{{\bf \thesection.\arabic{prop}.}\ \ \bf#1:}}{\smp}
\newenvironment{enonce}[1]{\pagebreak[3]\refstepcounter{prop}\mmpn{{\bf  \thesection.\arabic{prop}.\ #1.}}\begin{it} }{\end{it}\smp}
\def\thesection{\arabic{section}}

\newcommand{\result}[1]{\begin{enonce}{#1}}
\newcommand{\fresult}{\end{enonce}}
\def\dsp{\displaystyle}
\def\abeldeux{\mathcal{E}_{2}^\sharp}
\def\apdeux{\mathcal{A}_{\geq 2}}
\def\aprangdeux{\mathcal{A}_{=2}}
\newcommand{\lproj}[1]{\limproj{\un< Q\leq #1}D\big(N_#1(Q)/Q\big)}
\def\ls#1#2{{\hspace{.3ex}^{#1}\hspace{-.1ex}#2}}  
\begin{document}
\centerline{\large\bf Gluing endo-permutation modules}\vspace{.5cm}
\centerline{\bf Serge Bouc}
\vspace{1cm}\par

\begin{footnotesize}
{\bf Abstract :} In this paper, I show that if $p$ is an odd prime, and if $P$ is a finite $p$-group, then there exists an exact sequence of abelian groups
$$0\to T(P)\to D(P)\to\lproj{P}\to H^1\big(\apdeux(P),\Z\big)^{(P)}\mvirg$$
where $D(P)$ is the Dade group of $P$ and $T(P)$ is the subgroup of endo-trivial modules. \par
Here $\lproj{P}$ is the group of sequences of compatible elements in the Dade groups $D\big(N_P(Q)/Q\big)$ for non trivial subgroups $Q$ of $P$. The poset $\apdeux(P)$ is the set of elementary abelian subgroups of rank at least 2 of $P$, ordered by inclusion. The group $H^1\big(\apdeux(P),\Z\big)^{(P)}$ is the subgroup of $H^1\big(\apdeux(P),\Z\big)$ consisting of classes of $P$-invariant 1-cocycles.\par
A key result to prove that the above sequence is exact is a characterization of elements of $2D(P)$ by sequences of integers, indexed by sections $(T,S)$ of $P$ such that $T/S\cong (\Z/p\Z)^2$, fulfilling certain conditions associated to subquotients of $P$ which are either elementary abelian of rank~3, or extraspecial of order $p^3$ and exponent $p$.
\\ \vspace{1ex}\par
\vspace{-.5ex}\par
\end{footnotesize}
\section{Introduction}
The classification of all endo-permutation modules for finite $p$-groups has been completed recently, thanks to the work of several authors (see in particular \cite{alperin}, \cite{cath2},  \cite{cath3}, \cite{boma}, \cite{dadegroup}). This paper addresses the question of {\em gluing arbitrary endo-permutation modules}, and it is intended to be a complement to our previous joint work with Jacques Th\'evenaz~(\cite{BoTh2}), where the case of torsion endo-permutation modules was handled.\par
The gluing problem is the following~: let $p$ be an odd prime, let $P$ be a finite $p$-group, and let $k$ be a field of characteristic~$p$. If $v$ is an element of the Dade group $D(P)$ of endo-permutation $kP$-modules, and if $Q$ is a non trivial subgroup of $P$, denote by $v_Q$ the image of $v$ by the deflation-restriction map $\Defres_{N_P(Q)/Q}^P$. Then the $v_Q$'s are subject to some obvious compatibility conditions. Conversely, if $Q$ is a non-trivial subgroup of $P$, let $u_Q$ be an element of the Dade group $D_k\big(N_P(Q)/Q\big)$, and assume that these compatibility conditions between the $u_Q$'s are fulfilled. Is there an element $u\in D(P)$ such that for any non trivial subgroup $Q$ of $P$
$$\Defres_{N_P(Q)/Q}^{P}(u) = u_Q\;\;\;?$$ 
Such an element $u$ is called a {\em solution to the gluing problem} for the {\em gluing data} $(u_Q)_{1<Q\leq P}$.\par
When $P$ is abelian, the gluing problem was completely solved by Puig~\cite{puigcorr} (see also Lemma~\ref{sigma} below), and he used the result to construct suitable stable equivalences between blocks.\par
The main result of the present paper is that if $p$ is an odd prime, and if $P$ is a finite $p$-group, then there exists an exact sequence of abelian groups
$$0\longrightarrow T(P)\longrightarrow D(P)\longrightarrow\lproj{P}\stackrel{h_P}{\longrightarrow} H^1\big(\apdeux(P),\Z\big)^{(P)}\mvirg$$
where $D(P)$ is the Dade group of $P$ and $T(P)$ is the subgroup of endo-trivial modules. Here $\lproj{P}$ is the group of gluing data for $P$, i.e. the group of sequences of compatible elements in the Dade groups $D\big(N_P(Q)/Q\big)$ for non trivial subgroups $Q$ of~$P$. The poset $\apdeux(P)$ is the set of elementary abelian subgroups of rank at least 2 of $P$, ordered by inclusion. The group $H^1\big(\apdeux(P),\Z\big)^{(P)}$ is the subgroup of $H^1\big(\apdeux(P),\Z\big)$ consisting of classes of $P$-invariant 1-cocycles.\par
The main consequence of this result is that if $H^1\big(\apdeux(P),\Z\big)=\zero$, then the gluing problem always has a solution. Unfortunately, the map $h_P$ is not surjective in general, so when $H^1\big(\apdeux(P),\Z\big)\neq\zero$, not much can be said at the time for the gluing problem. In Section~\ref{exemple Xp5}, the example of the extraspecial group of order $p^5$ and exponent~$p$ is described in details. In this case, the group $H^1\big(\apdeux(P),\Z\big)^{(P)}$ is a free group of rank~$p^4$, and the image of $h_P$ has finite index in this group. In particular it is non zero, and the gluing problem does not always have a solution.\par
It could be true in general that $h_P$ always has finite cokernel, and this would be enough to show that if $H^1\big(\apdeux(P),\Z\big)^{(P)}\neq \zero$, then the image of $h_P$ is non zero, hence that the gluing problem does not always have a solution~: it is known indeed that the group $H^1\big(\apdeux(P),\Z\big)^{(P)}$ is a free abelian group, since the poset $\apdeux(P)$ has the homotopy type of a wedge of spheres (see~\cite{BoTh4}).  
\masubsect{Notation}
Throughout this paper, the symbol $p$ denotes an odd prime number, and $P$ denotes a finite $p$-group. Inclusion of subgroups will be denoted by $\leq$, and strict inclusion by $<$. Inclusion up to $P$-conjugation will be denoted by $\leq_P$.\par
A {\em section} $(T,S)$ of $P$ is a pair of subgroups of $P$ with $S\normal T$. The factor group~$T/S$ is the corresponding {\em subquotient} of $P$. If $(T,S)$ is a section of $P$, then $N_P(T,S)$ denotes $N_P(T)\cap N_P(S)$. \par
A class $\mathcal Y$ of $p$-groups is said to be {\em closed under taking subquotients} if for any $Y\in \mathcal{Y}$ and any section $(T,S)$ of~$Y$, any group isomorphic to $T/S$ belongs to $\mathcal{Y}$. If $\mathcal{Y}$ is such a class, and $P$ is a finite $p$-group, let $\mathcal{Y}(P)$ be the set of sections $(T,S)$ of $P$ such that $T/S\in\mathcal{Y}$. \par
The symbol $X_{p^3}$ denotes an extraspecial $p$-group of order $p^3$ and exponent $p$. The symbol $\mathcal{X}_3$ denotes the class of $p$-groups which are either elementary abelian of rank at most~3, or isomorphic to $X_{p^3}$.
Thus, the symbol $\mathcal{X}_3(P)$ denotes the set of sections $(T,S)$ of $P$ such that $T/S$ is elementary abelian of rank at most 3, or isomorphic to $X_{p^3}$. Let moreover $\abeldeux(P)$ denote the set of sections $(T,S)$ of $P$ such that $T/S\cong (\Z/p\Z)^2$. \par
If $P$ is a finite $p$-group, and $k$ is a field of characteristic $p$, let $D(P)$ denote the Dade group of endo-permutation $kP$-modules. The field $k$ does not appear in this notation, because it turns out that $D(P)$ is independent of $k$, at least when $p$ is odd (see~\cite{dadegroup} Theorem~9.5 for details).\par
When $(T,S)$ is a section of $P$, there is a {\em deflation-restriction map} $\Defres_{T/S}^P:D(P)\to D(T/S)$, which is the group homomorphism obtained by composing {\em the restriction map} $\Res_T^P:D(P)\to D(T)$, followed by the {\em deflation map} $\Def_{T/S}^T:D(T)\to D(T/S)$.\par
Recall that if $X$ is a finite $P$-set, there is a corresponding element $\Omega_X$ of the Dade group of $P$, called the syzygy of the trivial module relative to $X$ (or the {\em $X$-relative syzygy} for short) ~: it is defined as the class of the kernel of the augmentation map $kX\to k$ when this does make sense, and by 0 otherwise (see e.g. \cite{omegarel} for details). When $X$ is the set $P$ itself, on which $P$ acts by multiplication, the corresponding element will be denoted by $\Omega_{P/\un}$ or $\Omega_P$. \par

\masubsect{Contents} This paper is organized as follows~: 
\begin{itemize}
\item In Section~\ref{statement}, I state the main theorem (Theorem~\ref{suite exacte}), and this requires in particular the definition of some objects and maps between them. 
\item Section~\ref{biset functors} recalls some notation on biset functors, forgetful functors between categories of biset functors, and corresponding adjoint functors.
\item Section~\ref{rang 2} is devoted to the main tool~(Theorem~\ref{imdefres}) used in the proof of Theorem~\ref{suite exacte}, namely a characterization by linear equations of the image of the group $2D(P)$ by the deflation-restriction maps to all subquotients $T/S$ of $P$ which are elementary abelian of rank~2. This characterization may be a result of independent interest. 
\item Section~\ref{preuve} exposes the proof of Theorem~\ref{suite exacte}.
\item Finally, Section~\ref{exemple Xp5} focuses on the example of the extraspecial $p$-group of order $p^5$ and exponent $p$~: the reason for choosing this particular group is twofold~: it is is one of the smallest $p$-groups $P$ for which $H^1\big(\apdeux(P),\Z\big)\neq \zero$, and moreover the Dade group of this $p$-group is rather well known, thanks to our joint work with Nadia Mazza (\cite{boma}).
\end{itemize}
\section{Statement of the main theorem}\label{statement}
\result{Notation} If $P$ is a finite $p$-group, then $\apdeux(P)$ denotes the poset of elementary abelian subgroups of $P$ of rank at least~2. Let $\aprangdeux(P)$ denote the set of elementary abelian subgroups of rank 2 of $P$.\fresult
Recall that if the $p$-rank of $P$ is at least equal to~3, then all the elementary abelian subgroup of $P$ of rank at least~3 are in the same connected component of $\apdeux(P)$. This component is called {\em the big component}. It is obviously invariant under $P$-conjugation. Each of the other connected components, if there are any, consists of a single maximal elementary abelian subgroup of rank~2.
\result{Notation}
Denote by $\limproj{\un<Q\leq P}D\big(N_P(Q)/Q\big)$ the set of sequences $(u_Q)_{\un<Q\leq P}$, indexed by non trivial subgroups of $P$, where $u_Q\in D\big(N_P(Q)/Q\big)$, such that~:
\begin{itemize}
\item If $x\in P$, then $\ls{x}{u}_Q=u_{\ls{x}{Q}}$.
\item If $Q\normal R$, then $\Defres_{N_P(Q,R)/R}^{N_P(Q)/Q}u_Q=\Res_{N_P(Q,R)/R}^{N_P(R)/R}u_R$.
\end{itemize}
Denote by $r_P: D(P)\to \lproj{P}$ the map sending $v\in D(P)$ to the sequence $(\Defres_{N_P(Q)/Q}^Pv)_{\un<Q\leq P}$. \par
If $E$ is an abelian $p$-group, denote by $\sigma_E$ the map $\limproj{\un<F\leq E}D(E/F)\to D(E)$ defined by
$$\sigma_E(u)=-\sum_{\un<F\leq E}\mu(\un,F)\Inf_{E/F}^Eu_F\mvirg$$
where $\mu$ is the M\"obius function of the poset of subgroups of $P$.
\fresult
It has been shown by Puig (\cite{puigdadetf} 2.1.2) that the kernel of~$r_P$ is equal to the group $T(P)$ of endo-trival modules. Moreover, when $E$ is an abelian group, the map $r_E$ is surjective (\cite{puigcorr} Proposition~3.6). More precisely~:
\result{Lemma}\label{sigma} Let $E$ be an abelian $p$-group. Then $\sigma_E$ is a section of $r_E$, i.e. $r_E\sigma_E$ is equal to the identity map of $\limproj{\un<F\leq E}D(E/F)$.
\fresult
\pf Let $\un<G\leq E$. Then
\begin{eqnarray*}
\Def_{E/G}^E\sigma_E(u)&=&-\sum_{\un<F\leq E}\mu(\un,F)\Def_{E/G}^E\Inf_{E/F}^Gu_F\\
&=&-\sum_{\un<F\leq E}\mu(\un,F)\Inf_{E/FG}^{E/G}\Def_{E/FG}^{E/F}u_F\\
&=&-\sum_{\un<F\leq E}\mu(\un,F)\Inf_{E/FG}^{E/G}u_{FG}\\
&=&-\sum_{G\leq R\leq E}\big(\sumb{\un<F\leq R}{FG=R}\mu(\un,F)\big)\Inf_{E/R}^{E/G}u_R\mpoint
\end{eqnarray*}
Now if $G<R$
$$\sumb{\un<F\leq R}{FG=R\rule{0ex}{1.5ex}}\mu(\un,F)=\sumb{\un\leq F\leq R}{FG=R\rule{0ex}{1.5ex}}\mu(\un,F)\mvirg$$ 
and this is equal to zero, by a classical combinatorial lemma, since $G\neq \un$. And if $G=R$
$$\sumb{\un<F\leq R}{FG=R\rule{0ex}{1.5ex}}\mu(\un,F)=-\mu(\un,\un)+\sum_{\un\leq F\leq R}\mu(\un,F)=-1\mpoint$$
Thus $\Def_{E/G}^E\sigma_E(u)=\Inf_{E/G}^{E/G}u_G=u_G$, as was to be shown.\findemo
\result{Lemma}\label{r_E surjective} Let $E$ be an elementary abelian group of rank at least~2. Then the map~$r_E$ is surjective, and its kernel is the free abelian group of rank one generated by~$\Omega_{E/\un}$.
\fresult
\pf The kernel of $r_E$ is the group $T(E)$ of endo-trivial modules. Since $E$ is elementary abelian, this group is free of rank one, generated by $\Omega_{E/\un}$, by Dade's Theorem (\cite{dade2}~\cite{dade2II}). The surjectivity of $r_E$ follows from Lemma~\ref{sigma}.\findemo
\masubsect{Restriction and conjugation} The following construction has been introduced in~\cite{BoTh2}, for the torsion subgroup of the Dade group, but it works as well for the whole Dade group~: let $P$ be a finite $p$-group, and $H$ be a subgroup of $P$. If $u\in\!\!\!\lproj{P}$, then the sequence $(v_Q)_{\un<Q\leq H}$ defined by
$$v_Q=\Res_{N_H(Q)/Q}^{N_G(Q)/Q}u_Q$$
is an element of $\lproj{H}$, denoted by $\Res_H^Pu$. The map $u\mapsto\Res_H^Pu$ is a linear map $\lproj{P}$ to $\lproj{Q}$. The following is the analogue of Lemma~2.4 of~\cite{BoTh2}~:
\result{Lemma}\label{res} Let $H$ be a subgroup of $P$. The diagram
$$\carre{D(P)}{r_P}{\limproj{\un<Q\leq  P}D\big(N_P(Q)/ Q\big)}{\Res_H^P}{\Res_H^P}
{D(H)}{r_H}{\limproj{\un<Q\leq  H}D\big(N_H(Q)/Q\big)}$$
is commutative.
\fresult
\pf This is straightforward.\findemo
Similarly, if $x\in P$, denote by $c_{x,H}:D(H)\to D(\ls{x}{H})$ the conjugation by $x$, sending $v$ to $\ls{x}{v}$. If $u\in\lproj{H}$, then the sequence $(v_R)_{\un<R\leq\ls{x}{H}}$ defined by $v_R=c_{x,R^x}({u}_{R^x})$ is an element of $\limproj{\un<R\leq\ls{x}{H}}D\big(N_{\ls{x}{H}}(R)/R\big)$, that will be denoted by~$\ls{x}{u}$. The assignment $u\mapsto \ls{x}{u}$ is a linear map from $\lproj{H}$ to $\limproj{\un<R\leq\ls{x}{H}}D\big(N_{\ls{x}{H}}(R)/R\big)$, also denoted by $c_{x,H}$. Then~:
\result{Lemma}\label{conj} Let $H$ be a subgroup of $P$, and let $E$ be an abelian subgroup of $P$. The following diagrams are commutative~:
$$\xymatrix{
D(H)\ar[r]^-{r_H}\ar[d]_{c_{x,H}}&*!U(0.4){\lproj{H}}\ar[d]^{c_{x,H}}\\
D(\ls{x}{H})\ar[r]^-{r_{\ls{x}{H}}}&*!U(0.2){\rule{0ex}{3ex}\limproj{\un<R\leq\ls{x}{H}}D\big(N_{\ls{x}{H}}(R)/R\big)}
}
\qquad\xymatrix{
D(E)\ar[d]_{c_{x,E}}&*!U(0.4){\limproj{\un<R\leq E}D(E/R)}\ar[d]^{c_{x,E}}\ar[l]_-{\sigma_E}\\
D(\ls{x}{E})&*!U(0.2){\rule{0ex}{3ex}\limproj{\un<R\leq\ls{x}{E}}D(E/R)}\ar[l]_-{\sigma_{\ls{x}{E}}}
}
$$
\fresult
\pf This is also straightforward.\findemo
\masubsect{Construction of a map} Let $E$ and $F$ be elements of $\apdeux(P)$ such that $E<F$. If $v\in\lproj{P}$, consider the element
$$d_{E,F}=\Res_E^F\sigma_F\Res_F^Pv-\sigma_E\Res_E^Pv$$
of $D(E)$. Then by Lemma~\ref{r_E surjective} and Lemma~\ref{sigma}
\begin{eqnarray*}
r_E(d_{E,F})&=&r_E\Res_E^F\sigma_F\Res_F^Pv-r_E\sigma_E\Res_E^Pv\\
&=&\Res_E^Fr_F\sigma_F\Res_F^Pv-r_E\sigma_E\Res_E^Pv\\
&=&\Res_E^F\Res_F^Pv-\Res_E^Pv=0\mpoint\\
\end{eqnarray*}
By Lemma~\ref{r_E surjective}, there exists a unique integer $w_{E,F}$ such that
$$d_{E,F}=w_{E,F}\cdot\Omega_{E/\un}\mpoint$$
If $x\in P$, then it is clear from Lemma~\ref{conj} that $\ls{x}{d}_{E,F}=d_{\ls{x}{E},\ls{x}{F}}$, and it follows that $w_{\ls{x}{E},\ls{x}{F}}=w_{E,F}$. Moreover, if $E,F,G\in\apdeux(P)$ with $E<F<G$, then 
$$d_{E,F}+\Res_{E}^Fd_{F,G}=d_{E,G}\mvirg$$
hence $w_{E,F}+w_{F,G}=w_{E,G}$. 
In other words the function sending the pair $(E,F)$ of elements of $\apdeux(P)$, with $E<F$, to $w_{E,F}$, is a $P$-invariant 1-cocycle on $\apdeux(P)$, with values in $\Z$~:
\result{Notation} Let $P$ be a finite $p$-group. A {\em $P$-invariant 1-cocycle} on $\apdeux(P)$, with values in $\Z$, is a function sending a pair $(E,F)$ of elements of $\apdeux(P)$, with $E<F$, to an integer $w_{E,F}$, with the following two properties~:
\begin{enumerate}
\item If $x\in P$ and $E<F$ in $\apdeux(P)$, then $w_{\ls{x}{E},\ls{x}{F}}=w_{E,F}$.
\item If $E<F<G$ in $\apdeux(P)$, then $w_{E,F}+w_{F,G}=w_{E,G}$.
\end{enumerate}
The set $\Big(Z^1\big(\apdeux(P)\big)\Big)^P$ of $P$-invariants 1-cocycles is a group for addition of functions.\par
Denote by $\Big(B^1\big(\apdeux(P)\big)\Big)^P$ the subgroup of $\Big(Z^1\big(\apdeux(P)\big)\Big)^P$ consisting of cocycles~$w$ for which there exists a $P$-invariant function $E\mapsto m_E$ from $\apdeux(P)$ to $\Z$ such that
$$\forall E<F \in \apdeux(P),\;w_{E,F}=m_F-m_E\mpoint$$
Denote by $H^1\big(\apdeux(P),\Z\big)^{(P)}$ the factor group $\Big(Z^1\big(\apdeux(P)\big)\Big)^P/\Big(B^1\big(\apdeux(P)\big)\Big)^P$
\fresult
\begin{rem}{Remark}\label{invariant} One can show that the group $\Big(B^1\big(\apdeux(P)\big)\Big)^P$ is also equal to the set of elements $w$ of $\Big(Z^1\big(\apdeux(P)\big)\Big)^P$ for which there exists a (not necessarily $P$-invariant) function $E\mapsto m_E$ such that $w_{E,F}=m_F-m_E$ for any $E<F$ in $\apdeux(P)$. This is because if $E<F$ in $\apdeux(P)$, then $E$ and $F$ are the ``big component", which is $P$-invariant. Since $w$ is $P$-invariant, it follows that the function
$$x\in P\mapsto m_{\ls{x}{E}}-m_E$$
does not depend on the choice of $E$, and that it is a group homomorphism from $P$ to~$\Z$ (i.e. an element of $H^1(P,\Z)$). There are no non zero such homomorphisms, so $m$ is actually $P$-invariant. 
\end{rem}
\begin{rem}{Remark} On the other hand, one can consider the ordinary first cohomology group $H^1\big(\apdeux(P),\Z\big)$ of $\apdeux(P)$ over $\Z$, which is defined similarly to $H^1\big(\apdeux(P),\Z\big)^{(P)}$, but forgetting all conditions of $P$-invariance. Then the group $P$ acts on $H^1\big(\apdeux(P),\Z\big)$, and it follows from Remark~\ref{invariant} that $H^1\big(\apdeux(P),\Z\big)^{(P)}$ is a subgroup of the group $H^1\big(\apdeux(P),\Z\big)^{P}$ of $P$-invariant elements in $H^1\big(\apdeux(P),\Z\big)$. It might happen however that this inclusion is proper~: an argument similar to the one used in Remark~\ref{invariant} yields an element in $H^2(P,\Z)$, and this group need not be zero.
\end{rem}
\result{Notation} Let $P$ be a finite $p$-group. Denote by 
$$h_P:\lproj{P}\to H^1\big(\apdeux(P),\Z\big)^{(P)}$$
the map sending $v\in\lproj{P}$ to the class of the 1-cocycle $w$ defined by the following equality, for $E<F$ in $\apdeux(P)$~:
\begin{equation}\label{defcocycle}w_{E,F}\cdot\Omega_{E/\un}=\Res_E^F\sigma_F\Res_F^Pv-\sigma_E\Res_E^Pv\mpoint
\end{equation}
\fresult
\result{Proposition}\label{compozero} Let $P$ be a finite $p$-group. Then $h_P$ is a group homomorphism, and the composition 
$$D(P)\stackrel{r_P}{\longrightarrow}\lproj{P}\stackrel{h_P}{\longrightarrow} H^1\big(\apdeux(P),\Z\big)^{(P)}$$
is equal to 0.
\fresult
\pf Clearly, the definition of $h_P$ implies that it is a group homomorphism. Observe next that for any $E\in\apdeux(P)$, since $r_E\sigma_E$ is the identity map, the image of the map $\sigma_Er_E-\Id_{D(E)}$ is contained in the kernel of $r_E$. By Lemma~\ref{r_E surjective}, it follows that there is a unique linear form $s_E$ on $D(E)$, with values in $\Z$, such that
$$\sigma_Er_E(u)=u+s_E(u)\cdot\Omega_{E/\un}\mvirg$$
for any $u\in D(E)$. By Lemma~\ref{conj}, this definition clearly implies that if $x\in P$, then $s_{\ls{x}{E}}(\ls{x}{u})=s_E(u)$, for any $u\in D(E)$.\par
Now if $E<F$ in $\apdeux(P)$, and if $v=r_P(t)$, for $t\in D(P)$, Equation~\ref{defcocycle} becomes
\begin{eqnarray*}
w_{E,F}\cdot\Omega_{E/\un}&=&\Res_E^F\sigma_F\Res_F^Pr_P(t)-\sigma_E\Res_E^Pr_P(t)\\
&=&\Res_E^F\sigma_Fr_F\Res_F^Pt-\sigma_Er_E\Res_E^Pt\;\;\;\hbox{(by Lemma~\ref{res})}\\
&=&\Res_E^F\big(\Res_F^Pt+s_F(\Res_F^Pt)\cdot\Omega_{F/\un}\big)-\big(\Res_E^Pt+s_E(\Res_E^Pt)\cdot\Omega_{E/\un}\big)\\
&=&\big(s_F(\Res_F^Pt)-s_E(\Res_E^Pt)\big)\cdot\Omega_{E/\un}\mpoint
\end{eqnarray*}
Setting $m_E=s_E(\Res_E^Pt)$, for $E\in\apdeux(P)$, yields
$$w_{E,F}=m_F-m_E\mvirg$$
hence $w\in \Big(B^1\big(\apdeux(P)\big)\Big)^P$ (the $P$-invariance of $m$ follows easily from the above remark, or from Remark~\ref{invariant}). Thus $h_Pr_P(u)=0$, as was to be shown.\findemo
The main theorem of this paper is the following~:
\result{Theorem}\label{suite exacte} Let $P$ be a finite $p$-group. Then the sequence of abelian groups
$$0\longrightarrow T(P)\longrightarrow D(P)\stackrel{r_P}{\longrightarrow} \lproj{P}\stackrel{h_P}{\longrightarrow}H^1\big(\apdeux(P),\Z\big)^{(P)}$$
is exact. 
\fresult
The key point in this theorem is to show that the kernel of $h_P$ is equal to the image of $r_P$. This will be done in two steps~: first take an element $u\in \Ker\;h_P$, and show that $2u\in r_P\big(2D(P)\big)$. This amounts to replacing $D$ by $2D$, which is easier to handle, since it is torsion free. Next, write $2u=r_P(2v)$, for some $v\in D(P)$. Then $u-r_P(v)$ is an element in $\limproj{\un<Q\leq P}D_t\big(N_P(Q)/Q\big)$, and it has been shown in \cite{BoTh2} that such a sequence of compatible torsion elements can always be glued (i.e. it always lies in $r_P\big(D(P)\big)$, though possibly not in $r_P\big(D_t(P)\big)$).
\pagebreak[2]
\section{Biset functors}\label{biset functors}
The main ingredient in the proof of Theorem~\ref{suite exacte} is the formalism of {\em biset functors}. A short exposition of the notation and main results on this subject can be found in Section~2 of~\cite{BoTh3}, Section~3 of~\cite{burnsideunits}, or Section~3 of~\cite{dadegroup}.\par
Recall in particular that if $(T,S)$ is a section of the group $P$, and if $M$ is a biset functor, then the set $P/S$ is a $(P,T/S)$-biset, and the corresponding {\em induction-inflation} morphism
$$M(P/S):M(T/S)\to M(P)$$
is denoted by $\Indinf_{T/S}^P$. Similarly, the set $S\dom P$ is a $(T/S,P)$-biset, and the corresponding {\em deflation-restriction} map
$$M(S\dom P):M(P)\to M(T/S)$$
is denoted by $\Defres_{T/S}^P$.\par
This notation was already used for the Dade group, and it is coherent~: it was shown more generally in~\cite{BoTh} that, if $P$ and $Q$ are finite $p$-groups, and $U$ is a finite $(Q,P)$-bisets, one can define a natural group homomorphism $D(U): D(P)\to D(Q)$. For an arbitrary prime number $p$, this construction does not yield a structure of biset functor for the whole Dade group, but only on the correspondence sending a $p$-group $P$ to the subgroup $D^\Omega(P)$ of $D(P)$ generated by all the relative syzygies $\Omega_X$ obtained for various $P$-sets~$X$. However, in the case $p=2$, it was shown in~\cite{dadegroup} that $D=D^\Omega$, so $D$ is indeed a biset functor in this case.\par
The following additional notation was introduced in~\cite{BoTh3}~: a class $\mathcal Y$ of $p$-groups is said to be {\em closed under taking subquotients} if for any $Y\in \mathcal{Y}$ and any section $(T,S)$ of~$Y$, the corresponding subquotient $T/S$ belongs to $\mathcal{Y}$. If $\mathcal{Y}$ is such a class, and $P$ is a finite $p$-group, denote by $\mathcal{Y}(P)$ the set of sections $(T,S)$ of $P$ such that $T/S\in\mathcal{Y}$. One can consider biset functors defined only on~$\mathcal Y$, with values in abelian groups.
Let ${\mathcal F}_{\mathcal Y}$ be the category of all such functors,
 and let $\mathcal F$ denote the category of functors defined on all finite $p$-groups. These categories are abelian categories. \par
The obvious forgetful functor
$${\mathcal O}_{\mathcal Y} : {\mathcal F} \to {\mathcal F}_{\mathcal
Y}$$
admits left and right adjoints
$${\mathcal L}_{\mathcal Y} : {\mathcal F}_{\mathcal Y} \to {\mathcal
F} \quad\text{and}\quad
{\mathcal R}_{\mathcal Y} : {\mathcal F}_{\mathcal Y} \to {\mathcal
F} \mvirg$$
whose evaluations can be computed as follows (cf. \cite{BoTh3}~Theorem~1.2 for details, in particular on the direct and inverse limits appearing in this statement)~:
\result{Theorem}
With the notation above, for any functor $M\in{\mathcal F}_{\mathcal
Y}$,  we have :
$$
{\mathcal L}_{\mathcal Y}M(P) \cong \limind{(T,S)\in{\mathcal Y}(P)} M
(T/S)
\qquad\text{and}\qquad
{\mathcal R}_{\mathcal Y}M(P) \cong \limproj{(T,S)\in{\mathcal Y}(P)}
M(T/S)\mpoint
$$
\fresult
Moreover (\cite{BoTh3}~Corollary~6.17), for any biset functor $M$ and any $p$-group $P$, the unit map
$$\eta_{M,P}:M(P)\to {\mathcal R}_{\mathcal Y}{\mathcal O}_{\mathcal Y}(P)=\limproj{(T,S)\in{\mathcal Y}(P)}M(T/S) $$
is given by
$$\eta_{M,P}(u)_{T,S}=\Defres_{T/S}^Pu\mvirg$$
for any $u\in M(P)$ and any $(T,S)\in\mathcal{Y}(P)$.
\section{Image in subquotients of rank two}\label{rang 2}
\result{Lemma}\label{limiteB*} Let $\mathcal{X}$ be a non-empty class of finite $p$-groups, closed under taking subquotients. Let $A$ be any abelian group, and denote by $\hat{B}$ the functor $\Hom_\Z\big(B(-),A\big)$, where $B$ is the Burnside functor. Then the unit map
$$\beta : \hat{B}\to \mathcal{R}_\mathcal{X}\mathcal{O}_\mathcal{X}\hat{B}$$
is an isomorphism.
\fresult
\pf Indeed if $P$ is a $p$-group, if $\varphi\in\Hom_\Z\big(B(P),A\big)$, if $(T,S)$ is a section of $P$ and $X$ is a subgroup such that $S\leq X\leq T$, then
$$(\Defres_{T/S}^P\varphi)\big((T/S)/(X/S)\big)=\varphi\big(\Indinf_{T/S}^P(T/S)/(X/S)\big)=\varphi(P/X)\mpoint$$
Suppose that $\varphi\in\Ker\;\beta$. Since the class $\mathcal{X}$ is non empty and closed under taking subquotients, then it contains the trivial group, and for any subgroup $X$ of $P$, the section $(X,X)$ is in $\mathcal{X}(P)$. It follows in particular that
$$0=(\Defres_{X/X}^P\varphi)\big((X/X)/(X/X)\big)=\varphi(P/X)\mpoint$$
Thus $\varphi=0$, and $\beta$ is injective.\par
Conversely, let $\psi=(\psi_{T,S})_{(T,S)\in\mathcal{X}(P)}$ be an element of $\limproj{\mathcal{X}(P)}\Hom_\Z\big(B(T/S),A\big)$. Equivalently, for each section $(T,S)\in\mathcal{X}(P)$ and each subgroup $X$ with $S\leq X\leq T$, we have an element $\psi_{T,S}\big((T/S)/(X/S)\big)$ of $A$, fulfilling the following two conditions~:
\begin{itemize}
\item[$i)$] If $x\in P$, then $\psi_{\ls{x}{T},\ls{x}{S}}\big((\ls{x}{T}/\ls{x}{S})/(\ls{x}{X}/\ls{x}{S})\big)=\psi_{T,S}\big((T/S)/(X/S)\big)$.
\item[$ii)$] If $(T,S)\in\mathcal{X}(P)$ and $(T',S')\in\mathcal{X}(P)$, and if $X$ is a subgroup of $P$ such that $S\leq S'\leq X\leq  T'\leq T$, then
$\psi_{T,S}\big((T/S)/(X/S)\big)=\psi_{T',S'}\big((T'/S')/(X/S')\big)$.
\end{itemize}
Condition $i)$ implies in particular that the element $\psi_{X,X}\big((X/X)/(X/X)\big)$ is constant on the conjugacy class of $X$ in $P$. Hence we can define an element $\varphi\in\Hom_\Z\big(B(P),A\big)$ by setting
$$\varphi(P/X)=\psi_{X,X}\big((X/X)/(X/X)\big)\mvirg$$
for any subgroup $X$ of $P$. Now if $(T,S)\in \mathcal{X}(P)$ and if $S\leq X\leq T$
\begin{eqnarray*}
(\Defres_{T/S}^P\varphi)\big((T/S)/(X/S)\big)&=&\varphi(P/X)=\psi_{X,X}\big((X/X)/(X/X)\big)\\
&=&(\Defres_{X/X}^{T/S}\psi_{T,S})\big((X/X)/(X/X)\big)\\
&=&\psi_{T,S}\big(\Indinf_{(T/S)/(X/X)}^{T/S}(X/X)/(X/X)\big)\\
&=&\psi_{T,S}\big((T/S)/(X/S)\big)\mpoint
\end{eqnarray*}
Thus $\Defres_{T/S}^P\varphi=\psi_{T,S}$ for any $(T,S)\in\mathcal{X}(P)$. Equivalently $\beta(\varphi)=\psi$, so $\beta$ is surjective, hence it is an isomorphism.\findemo
\result{Lemma}\label{lim2D} When $p$ is odd, the map
$$
\varepsilon=\prod_{(T,S)\in\mathcal{X}_3(P)} \Defres_{T/S}^P : \; 2D(P) \longrightarrow \limproj
{(T,S)\in\mathcal{X}_3(P)} 2D(T/S)
$$
is an isomorphism.
\fresult
\pf Consider the short exact sequence of $p$-biset functors
$$0\to 2D\to D\to\F_2D\to 0\mpoint$$
Applying the functor $\limproj{\mathcal{X}_3(P)}$ yields the bottom line of the following commutative diagram with exact lines
\begin{equation}\label{diagram1}
\xymatrix{
0\ar[r]&2D(P)\ar[r]\ar[d]^\varepsilon&D(P)\ar[r]\ar[d]^\delta& {\F}_2D(P)\ar[d]^\gamma\ar[r]&0\\
0\ar[r]& *!U(0.5){\limproj{\mathcal{X}_3(P)}2D}\ar[r]&*!U(0.5){\limproj{\mathcal{X}_3(P)}D}\ar[r]& *!U(0.5){\limproj{\mathcal{X}_3(P)} \F_2D}\mvirg
}
\end{equation}
where the map $\delta$ is an isomorphism, by Theorem~1.1 of~\cite{BoTh3}. Moreover, by Corollary~1.5 of~\cite{boya}, there is an exact sequence of $p$-biset functors
$$0\to B^\times\to \F_2B^*\to \F_2D^\Omega\to 0\mvirg$$
where $B^\times$ is the functor of units of the Burnside ring, which is isomorphic to the constant functor $\Gamma_{\F_2}$ for $p$ odd. Moreover $D^\Omega=D$ in this case. Applying the functor $\limproj{\mathcal{X}_3(P)}$ to this sequence yields the bottom line of the following commutative diagram with exact lines
\begin{equation}\label{diagram2}
\xymatrix{
0\ar[r]&{\F}_2\ar[r]\ar[d]^\alpha&{\F}_2B^*(P)\ar[r]\ar[d]^\beta& {\F}_2D(P)\ar[d]^\gamma\ar[r]&0\\
0\ar[r]& *!U(0.5){\limproj{\mathcal{X}_3(P)}\Gamma_{\F_2}}\ar[r]&*!U(0.5){\limproj{\mathcal{X}_3(P)}\F_2B^*}\ar[r]& *!U(0.5){\limproj{\mathcal{X}_3(P)} {\F}_2D}\mvirg
}
\end{equation}
Since $\F_2B^*$ is naturally isomorphic to $\Hom_\Z\big(B(-),\F_2\big)$, the map $\beta$ is an isomorphism, by Lemma~\ref{limiteB*}. Now the group $\limproj{\mathcal{X}_3(P)}\Gamma_{\F_2}$ is the set of sequences $(u_{T,S})_{(T,S)\in\mathcal{X}_3(P)}$ fulfilling the two following conditions~:
\begin{itemize}
\item[$i)$] If $(T,S)\in\mathcal{X}_3(P)$ and $x\in P$, then $u_{\ls{x}{T},\ls{x}{S}}=u_{T,S}$.
\item[$ii)$] If $(T,S)\in\mathcal{X}_3(P)$ and $(T',S')\in\mathcal{X}_3(P)$ are such that $S\leq S'\leq T'\leq T$, then $u_{T,S}=u_{T',S'}$.
\end{itemize}
Applying this for the case $T'=S'=S$, and next for the case $T=T'=S'$, it follows that $u_{T,S}=u_{S,S}=u_{T,T}$ for any $(T,S)\in\mathcal{X}_3(P)$. Thus $u_{T,T}=u_{S,S}$ if $T/S\in\mathcal{X}_3$. Since $\mathcal{X}_3$ contains the cyclic group of order $p$, and since $P$ is a $p$-group, it follows that $u_{T,T}=u_{\un,\un}$ for any subgroup $T$ of $P$. Hence $u_{T,S}$ is constant, and the map $\alpha$ is an isomorphism.\par
Now the Snake's Lemma, applied to Diagram~\ref{diagram2}, shows that the map $\gamma$ is injective. And another application of this Lemma to Diagram~\ref{diagram1} shows that the map $\varepsilon$ is an isomorphism.\findemo
Recall that if $E$ is an elementary abelian group of rank~2, then $2D(E)$ is free of rank one, generated by $2\Omega_{E/\un}$. Thus if $u\in 2D(P)$, and if $(T,S)\in\abeldeux(P)$, then $\Defres_{T/S}^Pu$ is a multiple of $2\Omega_{T/S}$. The following theorem characterizes the sequences of integers $(v_{T,S})_{(T,S)\in\abeldeux(P)}$ which can be obtained that way from an element of $2D(P)$~:
\result{Theorem}\label{imdefres} Let $P$ be a $p$-group (with $p>2$). The map
$$\mathcal{D}_P:2D(P)\to\prod_{(T,S)\in\abeldeux(P)}\Z$$
sending an element $u\in 2D(P)$ to the sequence $\mathcal{D}_P(u)_{T,S}$ of integers defined by
$$\Defres_{T/S}^P(u)=\mathcal{D}_P(u)_{T,S}\cdot 2\Omega_{T/S}\mvirg$$
is injective, and its image is equal to the set of sequences $(v_{T,S})_{(T,S)\in\abeldeux(P)}$ fulfilling the following conditions~:
\begin{enumerate}
\item If $(T,S)\in\abeldeux(P)$ and $x\in P$, then $v_{T,S}=v_{\ls{x}{T},\ls{x}{S}}$.
\item If $(T,S)$ and $(T',S)$ are in $\abeldeux(P)$, if $T'\leq N_P(T)$, and $TT'/S\cong (\Z/p\Z)^3$, then
$$v_{T,S}+\sum_{S<X<T}v_{TT',X}=v_{T',S}+\sum_{S<X<T'}v_{TT',X}\mpoint$$
\item If $(T,S)$ and $(T',S)$ are in $\abeldeux(P)$, if $T'\leq N_P(T)$, and $TT'/S\cong X_{p^3}$, then
$$v_{T,S}\equiv v_{T',S}\;\;\;({\rm mod.} \;p)\mpoint$$
\end{enumerate}
\fresult
\pf Let $u\in\Ker\;\mathcal{D}_P$. Then $\Defres_{T/S}^Pu=0$ for any $(T,S)\in \abeldeux(P)$. Moreover $\Defres_{T/S}^Pu\in 2D(\Z/p\Z)=\zero$ when $(T,S)$ is a section of $P$ with $T/S\cong \Z/p\Z$. The detection theorem of Carlson and Th\'evenaz (\cite{cath2}~Theorem~13.1) shows that $u=0$. Hence~$\mathcal{D}_P$ is injective.\par
Now Condition~1 of Theorem~\ref{imdefres} holds obviously for the elements of the image of~$\mathcal{D}_P$. To check that Conditions~2 and~3 also hold, suppose that $(T,S), (T',S)\in\abeldeux(P)$, with $T'\leq N_P(T)$ and $|TT':P|=p^3$, and observe that setting $R=TT'$, the diagram
$$\xymatrix{
D(P)\ar[r]^-{\mathcal{D}_P}\ar[d]_-{\Defres_{R/S}^P}&*!U(0.4){\mathop{\dsp{\prod}}_{(T",S")\in\abeldeux(P)}\limits\Z}\ar[d]^-{\pi_{R,S}}\\
D(R/S)\ar[r]^-{\mathcal{D}_P}&*!U(0.4){\mathop{\dsp{\prod}}_{(T"/S,S"/S)\in\abeldeux(R/S)}\limits\Z}
}$$
is commutative, where $\pi_{R,S}$ is the projection map obtained by identifying sections $(T",S")$ of $P$ such that $S\leq S"\leq T"\leq R$ with sections $(T"/S,S"/S)$ of $R/S$. This shows that it is enough to suppose that the group $P$ is either elementary abelian of rank~3, or isomorphic to $X_{p^3}$. These special cases are detailed below.\par 
To prove that conversely, Conditions~1,~2, and~3 characterize the image of $\mathcal{D}_P$, observe first that by Theorem~1.1 of~\cite{BoTh3}, the map
$$
\prod_{(T,S)} \Defres_{T/S}^P : \; D(P) \longrightarrow \limproj{(T,S)} D(T/S) $$
is an isomorphism, where $T/S$ runs through all sections of~$P$ which
are either elementary abelian $p$-groups of rank~$\leq3$ or
extraspecial groups of order~$p^3$ and exponent~$p$.\par
Suppose that Theorem~\ref{imdefres} is true when $P$ is elementary abelian of rank at most 3, or extraspecial of exponent $p$. Let $P$ be an arbitrary $p$-group, and consider a sequence $v=(v_{T,S})_{(T,S)\in\abeldeux(P)}$ fulfilling the conditions of Theorem~\ref{imdefres}. \par
If $(V,U)\in\mathcal{X}_3(P)$, then the correspondence $(T,S)\mapsto (T/U,S/U)$ is a one to one correspondence between the set of elements $(T,S)$ of $\abeldeux(P)$ such that $U\leq S\leq T\leq V$, and  $\abeldeux(V/U)$. Through this bijection, the sequence of integers $v_{T,S}$, for $U\leq S\leq T\leq V$, yields a sequence of integers fulfilling the conditions of Theorem~\ref{imdefres} for the group $V/U$, hence an element in the image of the map $\mathcal{D}_{V/U}$. In other words, there is a unique element $w_{V,U}\in 2D(V/U)$ such that
$$\Defres_{T/S}^{V/U}w_{V,U}=v_{T,S}\cdot 2\Omega_{T/S}\mvirg$$
for all $(T,S)\in\abeldeux(P)$ with $U\leq S\leq T\leq V$.\par
Now the uniqueness of $w_{V,U}$ shows that $w_{\ls{x}{V},\ls{x}{U}}=\ls{x}{w}_{V,U}$ for any $x\in P$, and that $\Defres_{V'/U'}^{V/U}w_{V,U}=w_{V',U'}$ whenever $(V,U)$ and $(V',U')$ are in $\mathcal{X}_3(P)$ and $U\leq U'\leq V'\leq V$. In other words, the sequence $(w_{V,U})_{(V,U)\in\mathcal{X}_3(P)}$ is an element of $\limproj{\mathcal{X}_3(P)}2D$. By Lemma~\ref{lim2D}, there exists an element $t\in 2D(P)$ such that
$$w_{V,U}=\Defres_{V/U}^Pt\mvirg$$
for any $(V,U)\in\mathcal{X}_3(P)$. Then obviously $\mathcal{D}_P(t)=v$, and $v$ lies in the image of $\mathcal{D}_P$, as was to be shown.\par
So the only thing left to check is that Theorem~\ref{imdefres} holds when $P$ is elementary abelian of rank at most~3, or isomorphic to $X_{p^3}$. This is a case by case verification, using the following lemma~:
\result{Lemma} \label{defres Omega_X}Let $P$ be a finite $p$-group, and $X$ be a finite $P$-set. 
\begin{enumerate}
\item If $T/S$ is a section of $P$, then 
$$\Defres_{T/S}^P\Omega_X=\Omega_{X^S}\mvirg$$
where $X^S$ denotes the set of fixed points of $S$ on $X$, viewed as a $T/S$-set. 
\item If moreover $(T,S)\in\abeldeux(P)$, then
$$\Defres_{T/S}^P(2\Omega_X)=\big(\sumb{S\leq V\leq T}{X^V\neq\emptyset}\mu(S,V)\big)2\Omega_{T/S}\mpoint$$
In other words $\mathcal{D}_P(2\Omega_X)_{T,S}=\dsp{\sumb{S\leq V\leq T}{X^V\neq\emptyset}\mu(S,V)}$.
\end{enumerate}
\fresult
\pf Assertion~1 follows from Section~4 of \cite{omegarel}. For Assertion~2, note that by Assertion~1 and Lemma~5.2.3 of~\cite{omegarel}, since $T/S$ is abelian,
$$\Defres_{T/S}^P2\Omega_X=\sumb{S\leq U\leq V\leq T}{X^V\neq\emptyset}\mu(U,V) \cdot 2\Omega_{T/U}\mvirg$$
and that $2\Omega_{T/U}=0$ in $D(T/S)$ unless $U=S$.
\findemo
Now there are four cases~:\medp
$\bullet$ If $|P|\leq p$, there is nothing to do, since the map $\mathcal{D}_P$ is an isomorphism $\zero\to\zero$.\medp
$\bullet$ If $P\cong (\Z/p\Z)^2$, then $2D(P)\cong\Z$, and $\abeldeux(P)=\{(P,\un)\}$. In this case, there is no condition on the image of $\mathcal{D}_P$, and $\mathcal{D}_P$ is an isomorphism $\Z\to\Z$. So Theorem~\ref{imdefres} holds in this case.\medp
$\bullet$ If $P\cong(\Z/p\Z)^3$, then $\abeldeux(P)$ consists of $p^2+p+1$ sections $(P,R)$, for $|R|=p$, and $p^2+p+1$ sections $(Q,\un)$, for $|Q|=p^2$. The group $2D(P)$ is a free abelian group, with basis
$$\{2\Omega_{P/\un}\}\sqcup\{2\Omega_{P/R}\mid |R|=p\}\mpoint$$
The following arrays gives the values of the sequence $v=\mathcal{D}_P(u)$ for the element $u$ in its first column on the left~:
$$\begin{array}{c||c|c}
&v_{P,R}&v_{Q,\un}\\
\hline\hline
2\Omega_{P/\un}&0&1\\
\hline
2\Omega_{P/R'}&\left\{\begin{array}{cl}0&\hbox{if}\;R'\neq R\\1&\hbox{if}\;R'= R\end{array}\right.&\left\{\begin{array}{cl}0&\hbox{if}\;R'< Q\\1&\hbox{if}\;R'\not< Q\end{array}\right.\\
\hline
\end{array}
$$
The image of the element $u=m_\un\cdot2\Omega_{P/\un}+\mathop{\sum}_{|R'|=p}\limits m_{R'}\cdot2\Omega_{P/R'}$ by the map~$\mathcal{D}_P$ is equal to the sequence $v=(v_{T,S})$, where
\begin{equation}\label{eq1}v_{P,R}=m_R\ressort{1cm}v_{Q,\un}=m_\un+\sum_{R\not< Q}m_R
\end{equation}
If $Q\neq Q'$ are subgroups of order $p^2$ of $P$, then $QQ'=P$, and
$$v_{Q,\un}+\sum_{\un<X<Q}v_{P,X}=m_\un+\sum_{|X|=p}m_X=v_{Q',\un}+\sum_{\un<X<Q'}v_{P,X}\mvirg$$
so Condition~2 of \ref{imdefres} holds for the sections $(Q,\un)$ and $(Q',\un)$ of $P$. Since $P$ is abelian, Conditions~1 and~3 of~\ref{imdefres} are obviously satisfied.\par
Conversely, suppose that Condition~2 hold for a sequence $v=(v_{T,S})_{(T,S)\in\abeldeux(P)}$.
This sequence is in the image of $\mathcal{D}_P$ if and only if there exist integers $m_\un$, $m_{R'}$, for $|R'|=p$, such that \ref{eq1} hold.\par
The first equation gives $m_R=v_{P,R}$, and then the second one gives
$$m_\un=v_{Q,\un}-\sum_{R\not<Q}v_{P,R}\mpoint$$
This is consistent if the right hand side does not depend on $Q$, i.e. if for any subgroups $Q\neq Q'$ of order $p^2$ of $P$
$$v_{Q,\un}-\sum_{R\not<Q}v_{P,R}=v_{Q',\un}-\sum_{R\not<Q'}v_{P,R}\mvirg$$
or equivalently
$$v_{Q,\un}+\sum_{R<Q}v_{P,R}=v_{Q',\un}+\sum_{R<Q'}v_{P,R}\mpoint$$
This is precisely Condition~2 of~\ref{imdefres} for the section $(Q,\un)$ and $(Q',\un)$, since $QQ'=P$ in this case. Thus Theorem~\ref{imdefres} holds for $P\cong(\Z/p\Z)^3$.\medp
$\bullet$ If $P\cong X_{p^3}$, then $\abeldeux(P)$ consists of the section $(P,Z)$, where $Z$ is the centre of~$P$, and of $p+1$ sections $(Q,\un)$, where $Q$ is a subgroup of index $p$ in $P$. The group $D(P)$ is equal to $D^\Omega(P)$, since $p\neq2$, so it is generated by the elements $\Omega_{P/\un}$, $\Omega_{P/X}$, for $|X|=p$, and $\Omega_{P/Q}$, for $|Q|=p^2$, which have order 2 in $D(P)$. Thus $2D(P)$ is generated by the elements $2\Omega_{P/\un}$ and $2\Omega_{P/X}$, for $|X|=p$. The following array gives the values of the sequence $v=\mathcal{D}_P(u)$ for the element $u$ in its first column on the left, where $R$ denotes a non central subgroup of order $p$ of $P$~:
$$\begin{array}{c||c|c}
&v_{P,Z}&v_{Q,\un}\\
\hline\hline
2\Omega_{P/\un}&0&1\\
\hline
2\Omega_{P/Z}&1&0\\
\hline
2\Omega_{P/R}&0&\left\{\begin{array}{cl}1&\hbox{if}\;R\not<Q\\1-p&\hbox{if}\;R<Q\end{array}\right.\\
\end{array}$$
The values in this table can be computed using Lemma~\ref{defres Omega_X}~: for example
$$\Res_Q^P2\Omega_{P/R}=\big(\sumb{\un\leq V\leq Q}{\rule{0cm}{1.5ex}V\leq_PR}\mu(\un,V)\big)2\Omega_{Q/\un}\mpoint$$
If $R\nleq Q$, then there is only one term in the summation, for $V=\un$, and $\mu(\un,V)=1$ in this case. And if $R\leq Q$, then there are $p$ additional terms, obtained for the $p$ distinct conjugates $V$ of $R$ in $P$, and $\mu(\un,V)=-1$ for each of them. This gives the value $1-p$ in this case.\par
Now if $u=m_{\un}\cdot 2\Omega_{P/\un}+m_Z\cdot 2\Omega_{P/Z}+\mathop{\sum}_{[R]}\limits m_{R}\cdot 2\Omega_{P/R}$ (where the brackets around $R$ mean that $R$ runs through a set of representatives of conjugacy classes of non central subgroups of order $p$ of~$P$), then the sequence $v=\mathcal{D}_P(u)$ is given by~:
\begin{equation}\label{eq6}
v_{P,Z}=m_Z\ressort{1cm}v_{Q,\un}=m_\un +\sum_{[R]\not<Q}m_{R}+(1-p)\sum_{[R]<Q}m_{R}\mvirg
\end{equation}
The second equation is equivalent to
\begin{equation}\label{second}
v_{Q,\un}=m_\un+\sum_{[R]}m_{R}-p\sum_{[R]<Q}m_{R}\mpoint
\end{equation}
It follows that $v_{Q,\un}\equiv v_{Q',\un}$ (mod. $p$), for any subgroups $Q$ and $Q'$ of order $p^2$ in~$P$. This shows that Condition~3 of \ref{imdefres} holds for the sections $(Q,\un)$ and $(Q',\un)$ of $P$. Condition~2 is obviously satisfied in this case, since $P$ has no subquotient isomorphic to $(\Z/p\Z)^3$.\par
Suppose now conversely that a sequence $v=(v_{T,S})_{(T,S)\in\abeldeux(P)}$ is given, and that Conditions 1 and 3 of~\ref{imdefres} hold. Then $v$ lies in the image of $\mathcal{D}_P$ if and only if there exist integers $m_\un$, $m_Z$, $m_{R}$ (invariant by $P$-conjugation), such that Equations~\ref{eq6} hold.\par
The first equation in~\ref{eq6} gives $m_Z=v_{P,Z}$, and the second one gives
$$m_\un-v_{Q,\un}+\sum_{[R]}m_{R}=p\sum_{[R]<Q}m_{R}\mpoint$$
All subgroups of order $p$ of $Q$ different from $Z$ are conjugate in $P$. Denoting by $R_Q$ one of them, this equation becomes
\begin{equation}\label{verif3}m_\un-v_{Q,\un}+\sum_{[R]}m_{R}=pm_{R_Q}\mpoint
\end{equation}
Summing this relation over $Q$ yields
$$(p+1)m_\un-\sum_{Q}v_{Q,\un}+(p+1)\sum_{[R]}m_{R}=p\sum_{[R]}m_{R}\mvirg$$
thus
$$\sum_{[R]}m_{R}=\sum_{Q}v_{Q,\un}-(p+1)m_\un\mpoint$$
Now~\ref{verif3} yields
$$pm_{R_Q}=\sum_{Q'\neq Q}v_{Q',\un}-pm_\un\mpoint$$
By Condition~3, the sum $\mathop{\sum}_{Q'\neq Q}\limits v_{Q',\un}$ is congruent to $pv_{Q,\un}$ modulo $p$, i.e. to~0. Since $Q=R_QZ$, this gives finally
$$m_R=\frac{1}{p}(\sum_{Q'\neq RZ}v_{Q',\un})-m_\un\mpoint$$
Conversely, if this holds for any $R$, then equation~\ref{second} holds~: indeed, in this case
$$\sum_{[R]}m_R=\sum_{Q}v_{Q,\un}-(p+1)m_\un\mvirg$$
thus
\begin{eqnarray*}
m_\un+\sum_{[R]}m_{R}-p\sum_{[R]<Q}m_{R}&=&m_\un+\sum_{Q}v_{Q,\un}-(p+1)m_\un-pm_{R_Q}\\
&=&\sum_{Q}v_{Q,\un}-pm_\un-(\sum_{Q'\neq Q}v_{Q',\un})+pm_\un\\
&=&v_{Q,\un}\mpoint
\end{eqnarray*}
Thus \ref{second} holds, and Theorem~\ref{imdefres} also, when $P=X_{p^3}$.\findemo
\section{Proof of Theorem~\ref{suite exacte}}\label{preuve}
Let $P$ be a finite $p$-group. Clearly $T(P)$ is the kernel of $r_P$, and $\Im\;r_P\leq\Ker\;h_P$, by Proposition~\ref{compozero}. So the only thing to show is that this inclusion is an equality.\par
Let $u\in\Ker\;h_P$. It means that there exists a $P$-invariant function $E\mapsto m_E$ from $\apdeux(P)$ to $\Z$ such that for any $E<F$ in $\apdeux(P)$
$$w_{E,F}=m_E-m_F\mvirg$$
where the integer $w_{E,F}$ is defined by the equality
$$w_{E,F}\cdot\Omega_{E/\un}=\Res_E^F\sigma_F\Res_F^Pu-\sigma_E\Res_E^Pu\mpoint$$
In other words
$$\Res_E^F(m_F\cdot\Omega_{F/\un}+\sigma_F\Res_F^Pu)=m_E\cdot\Omega_{E/\un}+\sigma_E\Res_E^Pu\mpoint$$
Set $w_E=m_E\cdot\Omega_{E/\un}+\sigma_E\Res_E^Pu$, for $E\in\apdeux(P)$. Then $\Res_E^Fw_F=w_E$, for any $E<F$ in $\apdeux(P)$, and $\ls{x}{(w_E)}=w_{\ls{x}{E}}$ for any $x\in P$ and $E\in\apdeux(P)$. Moreover , for any $E\in\apdeux(P)$
$$r_E(w_E)=r_E\sigma_E\Res_E^Pu=\Res_E^Pu\mvirg$$
since $r_E(\Omega_{E/\un})=0$, and since $\sigma_E$ is a section of $r_E$. It means that for any subgroup $Y\neq \un$ of $E$
$$\Def_{E/Y}^Ew_E=\Res_{E/Y}^{N_P(Y)/Y}u_Y\mpoint$$
If $(T,S)\in \abeldeux(P)$, define an integer $v_{T,S}$ by
\setlength{\arraycolsep}{0.5ex}
\begin{equation}\label{defvts}\left\{\begin{array}{rclc}\Res_{T/S}^{N_P(S)/S}(2u_S)&=&v_{T,S}\cdot 2\Omega_{T/S}& \;\hbox{if}\;\; S\neq \un\\2w_T&=&v_{T,\un}\cdot 2\Omega_{T/\un}&\; \hbox{if}\;\; S=\un\end{array}\right.
\end{equation}
This sequence of integers $(v_{T,S})_{(T,S)\in\abeldeux(P)}$ satisfies some of the conditions of Theorem~\ref{imdefres}. Indeed~:\par
$\bullet$ If $x\in P$ and $(T,S)\in\abeldeux(P)$, then $v_{\ls{x}{T},\ls{x}{S}}=v_{T,S}$~: this is because $\ls{x}{(w_E)}=w_{\ls{x}{E}}$ for any $E\in\apdeux(P)$, and because $\ls{x}{(u_Q)}=u_{\ls{x}{Q}}$ for any subgroup $Q\neq \un$ of $P$. Thus Condition~1 of Theorem~\ref{imdefres} holds.\par
$\bullet$ Suppose that $(T,S)$ and $(T',S)$ are elements of $\abeldeux(P)$ such that $T\leq N_P(T')$. There are two cases to consider~:
\begin{itemize}
\item[(a)] If $S\neq \un$, then for any section $(V,U)\in\abeldeux\big(N_P(S)/S\big)$
\begin{eqnarray*}
v_{V,U}\cdot\Omega_{V/U}&=&\Res_{V/U}^{N_P(U)/U}(2u_U)\\
&=&\Res_{V/U}^{N_P(S,U)/U}\Res_{N_P(S,U)/U}^{N_P(U)/U}(2u_U)\\
&=&\Res_{V/U}^{N_P(S,U)/U}\Defres_{N_P(S,U)/U}^{N_P(S)/S}(2u_S)\\
&=&\Defres_{V/U}^{N_P(S)/S}(2u_S)\mpoint\\
\end{eqnarray*}
It follows that the sequence $(v_{V,U})_{(V,U)\in\abeldeux\big(N_P(S)/S\big)}$ is equal to $\mathcal{D}_{N_P(S)/S}(2u_S)$, hence it is in the image of the map $\mathcal{D}_{N_P(S)/S}$. Thus if $TT'/S\cong(\Z/p\Z)^3$, then Condition~2 of Theorem~\ref{imdefres} holds for the sections $(T,S)$ and $(T',S)$ of $N_P(S)/S$. And if $TT'/S\cong X_{p^3}$, then Condition~3 of Theorem~\ref{imdefres} holds, for a similar reason.
\item[(b)] If $S=\un$, then set $F=TT'$. If $F\cong(\Z/p\Z)^3$, then consider a section $(V,U)\in\abeldeux(F)$. If $U=\un$, then
$$\Defres_{V/U}^F2w_F=\Res_V^F2w_F=2w_V=v_{V,\un}\cdot 2\Omega_{V/\un}\mpoint$$
And if $U\neq \un$, then
\begin{eqnarray*}
\Defres_{V/U}^F2w_F&=&\Res_{V/U}^{F/U}\Def_{F/U}^F2w_F\\
&=&\Res_{V/U}^{F/U}\Res_{F/U}^{N_P(U)/U}2u_U\\
&=&\Res_{V/U}^{N_P(U)/U}2u_U\\
&=&v_{V,U}\cdot2\Omega_{V/U}\mpoint
\end{eqnarray*}
It follows that the sequence $(v_{V,U})_{(V,U)\in\abeldeux(F)}$ is equal to $\mathcal{D}_F(2w_F)$. In particular, Condition~2 of Theorem~\ref{imdefres} is fulfilled for the sections $(T,\un)$ and $(T',\un)$ of $F$.
\end{itemize}
Hence the sequence $(v_{T,S})_{(T,S)\in\abeldeux(P)}$ fulfills all the conditions of Theorem~\ref{imdefres}, except possibly Condition~3 for sections $(T,\un)$ and $(T',\un)$ such that $T\leq N_P(T')$ and $TT'\cong X_{p^3}$. This situation is handled by the following lemma~:
\result{Lemma}\label{coton} Let $P$ be a finite $p$-group, and $(v_{T,S})_{(T,S)\in \abeldeux(P)}$ be a sequence of integers such that~:
\begin{enumerate}
\item If $x\in P$ and $(T,S)\in\abeldeux(P)$, then $v_{\ls{x}{T},\ls{x}{S}}=v_{T,S}$.
\item If $(T,S)$ and $(T',S)$ are in $\abeldeux(P)$, if $T\leq N_P(T')$ and if $TT'/S\cong(\Z/p\Z)^3$, then
$$v_{T,S}+\sum_{S<X<T}v_{TT',X}=v_{T',S}+\sum_{S<X<T'}v_{TT',X}\mpoint$$
\item If $S\neq \un$, if $(T,S)$ and $(T',S)$ are in $\abeldeux(P)$, if $T\leq N_P(T')$ and if $TT'/S\cong X_{p^3}$, then
$$v_{T,S}\equiv v_{T',S}\;\hbox{\rm (mod. $p$)}\mpoint$$
\end{enumerate}
Then~:
\begin{itemize}
\item[(i)] If $T,T'\in\aprangdeux(P)$, if $T$ and $T'$ are in the same connected component of $\apdeux(P)$, if $T\leq N_P(T')$ and $TT'\cong X_{p^3}$, then $v_{T,\un}\equiv v_{T',\un}\;\hbox{\rm (mod. $p$)}$.
\item[(ii)] There exists a sequence of integers $(y_T)_{T\in\aprangdeux(P)}$ such that
\begin{itemize}
\item[(a)] If $x\in P$ and $T\in\aprangdeux(P)$, then $y_{\ls{x}{T}}=y_T$.
\item[(b)] If $T,T'\in\aprangdeux(P)$, if $T\leq N_P(T')$ and $TT'\cong(\Z/p\Z)^3$, then $y_T=y_{T'}$.
\item[(c)] If $T,T'\in\aprangdeux(P)$, if $T\leq N_P(T')$ and $TT'\cong X_{p^3}$, then 
$$y_T+v_{T,\un}\equiv y_{T'}+v_{T',\un}\;\hbox{\rm (mod. $p$)}$$.
\end{itemize}
\end{itemize}
\fresult
\pf The proof of Assertion $(i)$ goes by induction on $|P|$, starting with the case where $P$ is cyclic, where there is nothing to prove. Assume then that Hypotheses~1), 2), and~3) imply Assertion~1, for any $p$-group of order strictly smaller than $|P|$. Let $T$ and $T'$ be elementary abelian subgroups of rank~2 of $P$, such that $T\leq N_P(T')$ and $TT'\cong X_{p^3}$. Set $X=TT'$, and denote by $Z$ the centre of $X$.\par
If there is a proper subgroup $Q$ of $P$ containing $X$, and such that $T$ and $T'$ are in the same connected component of $\apdeux(Q)$, then 
$v_{T,\un}\equiv v_{T',\un}\;\hbox{\rm (mod. $p$)}$, by induction, since Hypotheses~1), 2), and~3) obviously hold for $Q$ is they hold for $P$. It is the case in particular if $\apdeux(Q)$ is connected.\par
Suppose that there exists a subgroup $C$ of order $p$ in $C_P(X)$, not contained in $X$ (i.e. different from $Z$). Then the center $T"=C\times Z$ of the subgroup $Q=C\times X$ of~$P$ is not cyclic. Hence $\apdeux(Q)$ is connected, and $Q$ contains $T$ and $T'$. Thus I can suppose that $Q=P$, and then $T"$ is equal to the centre of $P$. It is elementary abelian of rank~2. Moreover $TT"\cong(\Z/p\Z)^3$, since $T$ and $T"$ are elementary abelian of rank 2 and centralize each other, and since $T\cap T"=T\cap X\cap T"=Z$. Hypothesis~2, applied to the sections $(T,\un)$ and $(T',\un)$ of $P$ yields
\begin{equation}\label{TT"}v_{T,\un}-v_{T",\un}=\sum_{\un<F<T"}v_{TT",F}-\sum_{\un<F<T}v_{TT",F}\mpoint
\end{equation}
Now $TT"\normal P$ since $|P:TT"|=p$, and $T\normal P$, since $T\normal X$ and $C\leq C_P(X)$. Hence $P$ acts by conjugation on the set of subgroups $F$ such that $\un<F<T$, and $F=Z$ is the unique fixed point under this action. Now Hypothesis~1 implies that
$$\sum_{\un<F<T}v_{TT",F}\equiv v_{TT",Z}\;\hbox{\rm (mod. $p$)}\mvirg$$
and Equation~\ref{TT"} yields
\begin{equation}\label{TT"1}v_{T,\un}-v_{T",\un}\equiv\sumb{\un<F<T"}{\rule{0ex}{1.5ex}F\neq Z}v_{TT",F}\;\hbox{\rm (mod. $p$)}\mpoint
\end{equation}
The same argument applies with $T'$ instead of $T$, so
\begin{equation}\label{T'T"1}v_{T',\un}-v_{T",\un}\equiv\sumb{\un<F<T"}{\rule{0ex}{1.5ex}F\neq Z}v_{T'T",F}\;\hbox{\rm (mod. $p$)}\mpoint
\end{equation}
Now if $\un<F<T"$ and $F\neq Z$, the group $P/F$ has order $p^3$ and exponent $p$ (since~$P$ has exponent $p$), and it is non abelian (since $F\not\leq [P,P]=Z$). Hence $P/F\cong X_{p^3}$. Since $P=(TT")(T'T")$, Hypothesis~3, applied to the sections $(TT",F)$ and $(T'T",F)$ yields $v_{TT",F}\equiv v_{T'T",F}\;\hbox{\rm (mod. $p$)}$. This shows that the right hand sides of~\ref{TT"1} and~\ref{T'T"1} are congruent modulo $p$. So are the left hand sides, thus $v_{T,\un}-v_{T",\un}\equiv v_{T',\un}-v_{T",\un}\;\hbox{\rm (mod. $p$)}$, and $v_{T,\un}\equiv v_{T',\un}\;\hbox{\rm (mod. $p$)}$.\par
Hence I can suppose that $Z$ is the only subgroup of order $p$ of $C_P(X)$. In particular, the centre of $P$ is cyclic, and $Z$ is the only subgroup of order $p$ in this centre. Moreover, since $T\neq T'$ and $T,T'$ are in the same connected component of $\apdeux(P)$, the groups $T$ and $T'$ are not maximal elementary abelian subgroups, thus $P$ has $p$-rank at least equal to~3, and $T$ and $T'$ are in the big component $\mathcal{C}$ of $\apdeux(P)$.\par
In this case, there is a normal subgroup $T_0$ of $P$ which is elementary abelian of rank~2, and $T_0\in\mathcal{C}$. Moreover $T_0>Z$.\par
If $T_0\not\leq X$, then $T_0\cap X=Z$. Then $|T_0X|=p^4$, and $|T_0X:X|=p$. Thus $T_0$ normalizes $X$.
Moreover, if $Y$ is a subgroup of index $p$ of $X$, then $Y>Z$, and $|T_0Y|=p^3$, thus $|T_0Y:Y|=p$, and $T_0$ normalizes $Y$. It follows that the image of $T_0$ in the group $\Out(X)$ of outer automorphisms of $X$, which is isomorphic to $GL_2(\F_p)$, is a $p$-subgroup stabilizing every line. So this image is trivial, and $T_0$ acts on $X$ by inner automorphisms. Let $t\in T_0-X$. Then there exists $y\in X$ such that $y^{-1}t\in C_P(X)$. In particular $y^{-1}t$ centralizes $y$, so $t$ centralizes $y$, and then $(y^{-1}t)^p=(y^{-1})^pt^p=1$. Moreover $y\neq t$, since $t\notin X$. Hence $y^{-1}t$ has order $p$. Since $Z$ is the only subgroup of order $p$ of $C_P(X)$, it follows that $y^{-1}t\in Z$, so $t\in X$. This contradiction shows that $T_0\leq X$. \par
Since the congruences $v_{T,\un}\equiv v_{T_0,\un}\;\hbox{\rm (mod. $p$)}$ and $v_{T',\un}\equiv v_{T_0,\un}\;\hbox{\rm (mod. $p$)}$ imply the congruence $v_{T,\un}\equiv v_{T',\un}\;\hbox{\rm (mod. $p$)}$, it is enough to suppose that $T_0=T$, thus $T\normal P$. Let~$F$ be an elementary abelian subgroup of rank~3 of $P$ containing $T'$~: such a subgroup exists, since $T'$ is not a maximal element of $\apdeux(P)$. Set $T"=C_F(T)$. Then $|F:T"|$ divide $p$, since $F/T"$ is a $p$-subgroup of $\Aut(T)\cong GL_2(\F_p)$. Moreover $F\not\leq C_P(T)$, since $F>T'$. Thus $|F:T"|=p$, and $T"\cong(\Z/p\Z)^2$. Moreover $T'\neq T"$, since $T'\not\leq C_P(T)$, thus $F=T'T"$.\par
Now $F$ centralizes $T'$, and normalizes $T$. Thus $F$ normalizes $TT'=X$. Moreover $F\cap X=T'$, since $T'\leq F\cap X$, and since $F$ and $X$ are distinct subgroups of order $p^3$ of~$P$, for $F$ is abelian and $X$ is not. Hence $|FX:F|=p$, so $FX$ normalizes $F$. Thus $X$ normalizes $F$, and $X$ also normalizes $C_P(T)$ since $T\normal P$. It follows that $X$ normalizes $F\cap C_P(T)=C_F(T)=T"$. Obviously $X$ also normalizes its subgroup $T'$.\par
Hypothesis 2, applied to the sections $(T',\un)$ and $(T",\un)$ of $P$, yields
\begin{equation}\label{sum1}v_{T',\un}-v_{T",\un}=\sum_{1<Y<T"}v_{T'T",Y}-\sum_{\un<Y<T'}v_{T'T",Y}\mpoint
\end{equation}
Since $X$ normalizes $T'$ and $T"$, and since any subgroup of order $p$ normalized by $X$ is centralized by $X$, Hypothesis~1 yields
$$\sum_{1<Y<T"}v_{T'T",Y}\equiv\sumb{1<Y<T"}{\rule{0ex}{1.5ex}Y\leq C_P(X)}v_{T'T",Y}\;\hbox{\rm (mod. $p$)}\mpoint$$
But $T"\cap C_P(X)=Z$, since $Z\leq T"\cap C_P(X)$, and since $Z$ is the only subgroup of order~$p$ of $C_P(X)$. Thus
\begin{equation}\label{cong1}\sum_{1<Y<T"}v_{T'T",Y}\equiv v_{T'T",Z}\;\hbox{\rm (mod. $p$)}\mpoint
\end{equation}
The same argument, applied with $T'$ instead of $T"$, since $T'\cap C_P(X)=Z$, yields
\begin{equation}\label{cong2}\sum_{1<Y<T'}v_{T'T",Y}\equiv v_{T'T",Z}\;\hbox{\rm (mod. $p$)}\mpoint
\end{equation}
Now it follows from \ref{sum1}, \ref{cong1} and~\ref{cong2} that $v_{T',\un}-v_{T",\un}\equiv 0\;\hbox{\rm (mod. $p$)}$, i.e. 
\begin{equation}\label{cong3}v_{T',\un}\equiv v_{T",\un}\;\hbox{\rm (mod. $p$)}\mpoint
\end{equation}
Now the group $TT"$ is also elementary abelian of rank 3~: indeed, the group $T"$ centralizes $T$, and $T"\cap T=Z$ since $T"\cap T\geq Z$ and $T"\neq T$ for $T\not\leq F\leq C_P(T')$. Then Hypothesis~2, for the sections $(T,\un)$ and $(T",\un)$, yields
\begin{equation}\label{sum2}v_{T,\un}-v_{T",\un}=\sum_{1<Y<T"}v_{TT",Y}-\sum_{\un<Y<T}v_{TT",Y}\mpoint
\end{equation}
The group $X$ normalizes $T$ and $T"$, and $T\cap C_P(X)=Z=T"\cap C_P(X)$. The same argument as above yields
\begin{equation}\label{cong4}v_{T,\un}\equiv v_{T",\un}\;\hbox{\rm (mod. $p$)}\mpoint
\end{equation}
Thus $v_{T,\un}\equiv v_{T',\un}\;\hbox{\rm (mod. $p$)}$, by~\ref{cong3} and~\ref{cong4}, and this completes the proof of Assertion~$(i)$.\par
For Assertion~2, there is nothing to do if $P$ has no normal subgroup $T_0\cong(\Z/p\Z)^2$, since then $P$ is cyclic, and $\apdeux(P)=\emptyset$. If $P$ is not cyclic, fix such a normal subgroup $T_0$ of $P$, and denote by $\mathcal{C}$ the connected component of $T_0$ in $\apdeux(P)$. Thus $\mathcal{C}$ is the big component if $P$ has $p$-rank at least~3, and $\mathcal{C}=\{T_0\}$ otherwise. Define the sequence $(y_T)_{T\in\aprangdeux(P)}$ by
$$y_T=\left\{\begin{array}{cc}0&{\rm if}\;T\in\mathcal{C}\\v_{T_0,\un}-v_{T,1}&{\rm otherwise}\end{array}\right.$$
This sequence obviously fulfills Condition~$(a)$ of Lemma~\ref{coton}, by Hypothesis~1, and since $\mathcal{C}$ is invariant by $P$-conjugation. Now it $T,T'\in \aprangdeux(P)$, if $T\leq N_P(T')$ and $TT'\cong (\Z/p\Z)^3$, it follows that $P$ has $p$-rank at least~3, that $\mathcal{C}$ is the big component, and that $T,T'\in\mathcal{C}$. Thus $y_T=y_{T'}=0$, so Condition~$(b)$ of Lemma~\ref{coton} holds. Finally, if $T,T'\in \aprangdeux(P)$, if $T\leq N_P(T')$ and $TT'\cong X_{p^3}$, then there are three cases~:
\begin{itemize}
\item if $T$ and $T'$ are in $\mathcal{C}$, then $y_T=y_{T'}=0$, thus $y_T+v_{T,\un}=v_{T,\un}$, and $y_{T'}+v_{T',\un}=v_{T',\un}$. But $v_{T,\un}\equiv v_{T',\un}\;\hbox{\rm (mod. $p$)}$ in this case, by Assertion~1. Thus Condition~$(c)$ holds in this case.
\item if $T\in\mathcal{C}$ and $T'\notin\mathcal{C}$, then $y_T+v_{T,\un}=v_{T,\un}$, and $y_{T'}+v_{T',\un}=v_{T_0,\un}$. But now $T$ and $T_0$ are both in $\mathcal{C}$, so $v_{T,\un}\equiv v_{T_0,\un}\;\hbox{\rm (mod. $p$)}$ in this case, by Assertion~1. Thus Condition~$(c)$ holds in this case also. The case $T\notin\mathcal{C}$ and $T'\in\mathcal{C}$ is similar.
\item if $T\notin\mathcal{C}$ and $T'\notin\mathcal{C}$, then $y_{T}+v_{T,\un}=v_{T_0,\un}=y_{T'}+v_{T',\un}$, so Condition~$(c)$ holds in this case also.
\end{itemize}
This completes the proof of Lemma~\ref{coton}.
\findemo
\noindent{\bf End of the proof of Theorem~\ref{suite exacte}~:} In the beginning of the proof of Theorem~\ref{suite exacte}, I started with an element $u\in\lproj{P}$ such that $h_P(u)=0$. From this data, in~\ref{defvts}, I built a sequence of integers $(v_{T,S})_{(T,S)\in\abeldeux(P)}$ fulfilling Hypothesis 1, 2 and~3 of Lemma~\ref{coton}. Let $(y_T)_{T\in\aprangdeux(P)}$ denote the sequence of integers provided by this lemma, and define a new sequence of integers $(v'_{T,S})_{(T,S)\in\abeldeux(P)}$ by
$$v'_{T,S}=\left\{\begin{array}{cc}v_{T,S}&{\rm if}\;S\neq \un\\y_T+v_{T,S}&{\rm if}\;S=\un\end{array}\right.\mpoint$$
Then this sequence fulfills Conditions~1, 2, and~3 of Theorem~\ref{imdefres}~: indeed, the new sequence is clearly invariant by conjugation, so Condition~1 is fulfilled. Conditions~2 and~3 for sections $(T,S)$ and $(T',S)$ with $S\neq \un$ are obviously fulfilled, since they are for the sequence $(v_{T,S})$, and since $v_{T,S}=v'_{T,S}$ when $S\neq\un$.\par
Now if $T,T'\in\aprangdeux(P)$, if $T\leq N_P(T')$ and $TT'\cong(\Z/p\Z)^3$, then $T$ and $T'$ are in the same connected component of $\apdeux(P)$, and $y_T=y_{T'}$. Thus
\begin{eqnarray*}
v'_{T,\un}+\sum_{\un<Y<T}v'_{TT',Y}&=&y_T+v_{T,\un}+\sum_{\un<Y<T}v_{TT',Y}\\
&=&y_{T'}+v_{T',\un}+\sum_{\un<Y<T'}v_{TT',Y}\\
&=&v'_{T',\un}+\sum_{\un<Y<T'}v'_{TT',Y}\mvirg
\end{eqnarray*}
so Condition~2 is fulfilled.\par
Finally if $T,T'\in\aprangdeux(P)$, if $T\leq N_P(T')$ and $TT'\cong X_{p^3}$, then
$$v'_{T,\un}=y_T+v_{T,\un}\equiv y_{T'}+v_{T',\un}\;\hbox{(mod. $p$)}\mvirg$$
hence Condition~3 is fulfilled, since $y_{T'}+v_{T',\un}=v'_{T',\un}$.\par
By Theorem~\ref{imdefres}, there exists $n\in D(P)$ such that $\mathcal{D}_P(2n)=(v'_{T,S})$. In other words, for any $(T,S)\in \abeldeux(P)$
$$\Defres_{T/S}^P(2n)=v'_{T,S}\cdot 2\Omega_{T/S}\mpoint$$
Thus if $S\neq \un$
$$\Defres_{T/S}^P(2n)=v_{T,S}\cdot 2\Omega_{T/S}=\Res_{T/S}^{N_P(S)/S}(2u_S)\mpoint$$
Set $t_S=\Defres_{N_P(S)/S}^P(2n)-2u_S$. Now for any $(V,U)\in\abeldeux\big(N_P(S)/S\big)$
\begin{eqnarray*}
\Defres_{V/U}^{N_P(S)/S}(t_S)&=&\Defres_{V/U}^P(2n)-\Defres_{V/U}^{N_P(S)/S}(2u_S)\\
&=&\Defres_{V/U}^P(2n)-\Res_{V/U}^{N_P(S,U)/U}\Defres_{N_P(S,U)/U}^{N_P(S)/S}(2u_S)\\
&=&\Defres_{V/U}^P(2n)-\Res_{V/U}^{N_P(S,U)/U}\Res_{N_P(S,U)/U}^{N_P(U)/U}(2u_U)\\
&=&\Defres_{V/U}^P(2n)-\Res_{V/U}^{N_P(U)/U}(2u_U)=0\mpoint
\end{eqnarray*}
It follows that $t_S$ is a torsion element of $D\big(N_P(S)/S\big)$, which is also in $2D\big(N_P(S)/S\big)$. Since the latter is torsion free, it follows that $t_S=0$, i.e. that $2u_S=\Defres_{N_P(S)/S}^P(2n)$, for any $S\neq \un$. Equivalently $2u=r_P(2n)$, or $2\big(u-r_P(n)\big)=0$.\par
Now $u-r_P(n)$ is an element of $\limproj{\un<Q\leq P}D_{tors}\big(N_P(Q)/Q\big)$. By Proposition~5.5 of~\cite{BoTh2}, there exists an element $m\in D(P)$ such that $r_P(m)=u-r_P(n)$. It follows that $u=r_P(m+n)$, as was to be shown. This completes the proof of Theorem~\ref{suite exacte}.\findemo
\section{Example~: the group $X_{p^5}$}\label{exemple Xp5}
Let $P$ be an extraspecial group of order $p^5$ and exponent $p$. The centre $Z$ of $P$ is cyclic of order $p$, and it is equal to the Frattini subgroup of $P$. The commutator $P\times P\to Z$ induces a non degenerate symplectic $\F_p$-valued scalar product on the factor group $E=P/Z\cong (\F_p)^4$, and the map $Q\mapsto Q/Z$ is a poset isomorphism from the poset of elementary abelian subgroups of $P$ strictly containing $Z$ to the poset $\mathcal{E}$ of non zero totally isotropic subspaces of $E$. There are $e=\frac{\dsp p^4-1}{\dsp p-1}$ isotropic lines in $E$, and the same number of totally isotropic 2-dimensional subspaces. It follows that $|\mathcal{E}|$ is equal to $2e$.\par
There is a commutative diagram
$$\xymatrix{
&&D(P/Z)\ar[r]_-\cong^-d&*!U(0.4){\limproj{Q\geq Z}D(P/Q)}\\
0\ar[r]&T(P)\ar[r]&D(P)\ar[r]^-{r_P}\ar[u]^{\Def_{P/Z}^P}&*!U(0.2){\rule{0ex}{3ex}\lproj{P}}\ar[r]^{h_P}\ar[u]_{\pi}&H^1\big(\apdeux(P),\Z\big)^{(P)}\\
}
$$
In this diagram, the group $\limproj{Q\geq Z}D(P/Q)$ is the group of sequences $(u_Q)_{Z\leq Q\leq P}$, where $u_Q\in D(P/Q)$ (note that $Q\normal P$ if $Q\geq Z$), such that
$$\forall R\geq Q\geq Z,\;\Def_{P/R}^{P/Q}u_Q=u_R\mvirg$$
and the map $\pi$ is the projection map on the components $Q\geq Z$.
 The map $d$ is the product of the deflation maps $\Def_{P/Q}^{P/Z}$, for $Q\geq Z$. It is an isomorphism, since the sequence $(u_Q)_{Q\geq Z}$, where $u_Q\in D(P/Q)$, is in the group $\limproj{Q\geq Z}D(P/Q)$ if and only if $u_Q=\Def_{P/Q}^{P/Z}u_Z$ for any $Q\geq Z$.\par
The kernel of the map $\Def_{P/Z}^P$ is the set of faithful elements of $D(P)$, and it was denoted by $\partial D(P)$ in \cite{boma}. It was shown in this paper (Theorem~9.1) that 
$$\partial D(P)\cong \Z^{2e}\oplus \Z/2\Z\mpoint$$ 
The kernel of $\pi$ consists of the sequences $(u_Q)_{\un<Q\leq P}$ in $\lproj{P}$ for which $u_Q=0$ if $Q\geq Z$. It was shown in~\cite{boma} that this is also the group $\limproj{{\scriptstyle \un<Q\leq P}\atop{\scriptstyle Q\cap Z=\un\rule{0ex}{1.5ex}}}\partial D\big(N_P(Q)/Q\big)$.\par
All these facts show that there is an exact sequence
\begin{equation}\label{xp5}0\longrightarrow T(P)\longrightarrow \partial D(P)\stackrel{r_P}{\longrightarrow} \limproj{{\scriptstyle \un<Q\leq P}\atop{\scriptstyle Q\cap Z=\un\rule{0ex}{1.5ex}}}\partial D\big(N_P(Q)/Q\big)\stackrel{h_P}{\longrightarrow} H^1\big(\apdeux(P),\Z\big)^{(P)}\mvirg
\end{equation}
where $r_P$ and $h_P$ are the restrictions of the previously defined maps with the same names to the corresponding subgroups.\par
If $Q$ is a subgroup of $P$ such that $Q\cap Z=\un$, then $Q$ is elementary abelian of rank at most 2 (see~\cite{boma} for details). If $Q$ has order $p$, then $N_P(Q)/Q\cong X_{p^3}$, thus $\partial D\big(N_P(Q)/Q\big)\cong \Z^{p+1}\oplus \Z/2\Z$ (\cite{boma} Theorem~9.1 or Section~11). If $Q$ has order $p^2$, then $N_P(Q)/Q\cong \Z/p\Z$, thus $\partial D\big(N_P(Q)/Q\big)\cong\Z/2\Z$. \par
It follows easily that the group $\limproj{{\scriptstyle \un<Q\leq P}\atop{\scriptstyle Q\cap Z=\un\rule{0ex}{1.5ex}}}\partial D\big(N_P(Q)/Q\big)$ has free rank at least equal to $e(p+1)$, since it contains the group $\mathop{\oplusb{|Q|=p}{Q\neq Z}\limits }2\partial D\big(N_P(Q)/Q\big)$.\par
Now the group $T(P)$ is free of rank one, generated by $\Omega_{P/\un}$, by Corollary~1.3 of~\cite{cath2}, and the group $H^1\big(\apdeux(P),\Z\big)^{(P)}$ is isomorphic to $H^1(\mathcal{E},\Z)$. An easy computation, using e.g. Section~6 of~\cite{boma}, shows that this group is free of rank~$p^4$.\par
Now the free rank of the image of $h_P$ in the exact sequence~\ref{xp5} is at least equal to
$$1-2e+e(p+1)=1+e(p-1)=p^4\mvirg$$
and since this is equal to the free rank of $H^1\big(\apdeux(P),\Z\big)^{(P)}$, it follows that the free rank of the image of $h_P$ is actually equal to $p^4$, and that the free rank of $\limproj{{\scriptstyle \un<Q\leq P}\atop{\scriptstyle Q\cap Z=\un\rule{0ex}{1.5ex}}}\partial D\big(N_P(Q)/Q\big)$ is equal to $e(p+1)$.\par
Moreover the map $h_P$ has finite cokernel, and this shows that in this case, the gluing problem does not always have a solution.
\begin{rem}{Remark} In this case, a precise description of the map $h_P$ shows that its cokernel is a non trivial finite $p$-group.
\end{rem}
\begin{rem}{Aknowledgements} I wish to thank Nadia Mazza and Jacques Th\'evenaz for careful reading of an early version of this paper, and for many suggestions, comments, and stimulating discussions about it.
\end{rem}

\par\noindent
\begin{flushleft}
Serge Bouc, CNRS - UMR 6140 - LAMFA, Universit\'e de Picardie - Jules Verne, 33, rue St Leu, F-80039 Amiens Cedex~1, France. 
\par\noindent {email~: {\tt serge.bouc@u-picardie.fr}}
\end{flushleft}
\end{document}